\documentclass[square]{ws-procs975x65}

\begin{document}

\title{Equivalence of Sobolev Inequalities and Lieb--Thirring Inequalities\footnote{\copyright\, 2009 by the authors. This paper may be reproduced, in its entirety, for non-commercial purposes.\\
This work was partially supported by U.S. NSF grants PHY-0652854 (R.L.F. and E.H.L.) and PHY-0652356 (R.S.).}}

\author{Rupert L. Frank}
\address{Department of Mathematics,
Princeton University, Washington Road, Princeton, NJ 08544, USA \\ Email: rlfrank@math.princeton.edu}

\author{Elliott H. Lieb}
\address{Departments of Mathematics and Physics,
Princeton University,
       P.~O.~Box 708, Princeton, NJ 08542, USA \\ Email: lieb@princeton.edu}

\author{Robert Seiringer}
\address{Department of Physics, Princeton University,
P.~O.~Box 708,
       Princeton, NJ 08542, USA\\ Email: rseiring@princeton.edu}

\begin{abstract}
We show that, under very general definitions of a kinetic energy operator $T$, the Lieb--Thirring inequalities for sums of eigenvalues of $T-V$ can be derived from the 
Sobolev inequality appropriate to that choice of $T$. 

\end{abstract}

\keywords{Schr\"odinger operator, Sobolev inequality, bound states, stability of matter}

\bodymatter


\section{Introduction}

The Sobolev and Lieb--Thirring (LT) inequalities seem to be very
different. The former is a kind of uncertainty principle, which,
effectively, states how large a negative potential $-V$ must be for
the Schr\"odinger operator $H=-\Delta -V$ to have a bound state. The
latter, which was originally introduced to prove the stability of
matter \cite{LiSe}, estimates the sum of all the negative eigenvalues
of $H$ and, apparently, is stronger.  Our goal here is to summarize
some recent work that says, surprisingly, that the latter can,
nevertheless, be derived from the former.

This theme extends to other, more complicated operators than $T=
-\Delta$ . For example, to the barely positive ``Hardy'' operator
$T=-\Delta - (d-2)^2/(4|x|^2)$ for dimensions $d\geq 3$. Another is
the relativistic energy $T=\sqrt{-\Delta + m^2} -m$, or just
$(-\Delta)^s$.  Still another is the inclusion of a magnetic vector
potential $A$ and $T=-(\nabla +iA(x))^2$. In all cases there is a
Sobolev type inequality $(\psi, T \, \psi) \geq C\, \Vert \psi
\Vert^2$, where $\Vert \cdot\Vert $ is a suitable norm. 

The LT inequalities \cite{LiTh} are of the form $\sum_j
|\lambda_j|^\gamma \leq C'\, (\Vert V \Vert')^p$, where the
$\lambda_j$ are the negative eigenvalues of $H=T-V$, $\Vert V \Vert'$
is another norm, and the allowed range of exponents $\gamma$ and $p$
depends on the dimension $d$. A bound in the limiting case $\gamma
=0$, i.e., a bound on the number of bound states of $-\Delta-V$, which
is valid for $d\geq 3$, is due to Cwikel \cite{Cw}, Lieb \cite{Li} and
Rozenblum \cite{Ro} and is called the CLR bound.

We shall explain, in very general terms, how LT and CLR inequalities
can be derived from Sobolev inequalities (the converse being almost
trivial).  We shall also give several examples for the purpose of
clarification.  One of the most physically relevant of these is to the
relativistic $T$, which in connection with earlier work
\cite{Lieb-Yau} yields a proof of the stability of relativistic matter
in arbitrary magnetic fields all the way up to the critical value of
the allowed nuclear charge $Z= 2/\pi\alpha$, where $\alpha $ is the
fine structure constant.


\section{Main results}

\subsection{The setup}
We start with an abstract setting, just to show how general the equivalence of Sobolev and 
LT is. Much of this can be skipped for practical applications. 

Let $X$ be a sigma-finite measure space. We consider the measure on $X$ as fixed and denote integration with respect to this measure by $dx$. By $L^p(X)$ for $1\leq p\leq \infty$ we denote the usual $L^p$ space with respect to this measure. Moreover, if $w$ is a non-negative function on $X$ we write $L^p(X,w)$ for the $L^p$ space with respect to the measure $w(x) dx$.

Let $t$ be a non-negative quadratic form, with domain $\text{dom}\, t$, which is closed in the Hilbert space $L^2(X)$ and let $T$ be the corresponding self-adjoint operator.

Throughout this paper we work under the following assumption which depends on a parameter $1<\kappa<\infty$.

\begin{assumption}[Generalized Beurling-Deny conditions]\label{ass:bd}
$\phantom{x}$
\vspace{-.2truecm}
\begin{enumerate}
 \item 
if $u,v\in\mathrm{dom}\, t$ are real-valued then $t[u+iv]=t[u]+t[v]$,
\item
if $u\in\mathrm{dom}\, t$ is real-valued then $|u|\in\mathrm{dom}\, t$ and $t[|u|]\leq t[u]$.
\item
there is a measurable, a.e. positive function $\omega$ such that if $u\in\mathrm{dom}\, t$ is non-negative then $\min(u,\omega)\in\mathrm{dom}\, t$ and $t[\min(u,\omega)]\leq t[u]$. Moreover, there is a form core $\mathcal Q$ of $t$ such that $\omega^{-1}\mathcal Q$ is dense in $L^2(X, \omega^{2\kappa/(\kappa-1)})$.
\end{enumerate}
\end{assumption}

For $\omega\equiv 1$, these are the usual Beurling--Deny conditions; see, e.g., \cite[Sec. 1.3]{D}. We note that in this case the assumption is independent of the value of $\kappa$. In our applications below it will be important to allow for $\omega\not\equiv 1$. In those examples, $X$ will have a differentiable structure and the density assumption will be satisfied for all $\kappa$ because $\omega$ is sufficiently smooth. We refer to Section \ref{sec:examples} for those examples.


\subsection{Main results}

Our first result concerns upper bounds on the number of negative eigenvalues of Schr\"odinger operators $T-V$ in terms of integrals of the potential $-V$. We denote by $N(-\tau,T-V)$ the number of eigenvalues less than $-\tau\leq 0$ of the operator $T-V$, taking multiplicities into account, and we abbreviate $N(T-V) :=N(0,T-V)$. We shall prove

\begin{theorem}[Equivalence of Sobolev and CLR Inequalities]\label{clrweighted}
 Under Assumption \ref{ass:bd} for some $\kappa>1$ the following are equivalent:
\begin{enumerate}
 \item\label{it:sobolevweighted} $T$ satisfies a Sobolev inequality with exponent $q=2\kappa/(\kappa-1)$, that is, there is a constant $S>0$ such that for all $u\in\mathrm{dom}\, t$,
\begin{equation}\label{eq:sobolevweighted}
 t[u] \geq S \left( \int_X |u|^q \,dx \right)^{2/q} \,.
\end{equation}
\item\label{it:clrweighted} $T$ satisfies a CLR inequality with exponent $\kappa$, that is, there is a constant $L>0$ such that for all $0\leq V\in L^\kappa(X)$,
\begin{equation}\label{eq:clrweighted}
 N(0,T-V) \leq L \int_X V^{\kappa} \,dx \,.
\end{equation}
\end{enumerate}
The respective constants are related according to
\begin{equation}\label{eq:CLRconst}
S^{-\kappa} \leq L  \leq e^{\kappa-1} S^{-\kappa} \,.
\end{equation}
\end{theorem}

We emphasize that the statement of the theorem does not depend on $\omega$ in Assumption~\ref{ass:bd}. Only its existence and not its form is relevant. 

The implication $(\ref{it:clrweighted}\Rightarrow\ref{it:sobolevweighted})$ is a simple consequence of H\"older's inequality and the variational characterization of the lowest eigenvalue and is valid \emph{without} Assumption~\ref{ass:bd}. The converse is much deeper. 

\begin{remark}
Below we shall sketch two proofs of $(\ref{it:sobolevweighted}\Rightarrow\ref{it:clrweighted})$  which are abstract versions of proofs by Lieb \cite{Li} and by Li and Yau \cite{LYa}; see also \cite{LeSo} for the case $\omega\equiv 1$.
The latter method gives the bound $L  \leq e^{\kappa-1} S^{-\kappa}$ stated in \eqref{eq:CLRconst}. The method of \cite{Li} proceeds via the bound
\begin{equation*}
 \|\exp(-sT)\|_{1\to \infty} \leq K s^{-\kappa} 
\end{equation*}
which follows from the Sobolev inequality (\ref{eq:sobolevweighted}) with $K\leq (\kappa/S)^\kappa$. (For simplicity we consider the case $\omega\equiv 1$ here.) The method then gives the bound 
\begin{equation*}
L  \leq  \frac{K}{\kappa(\kappa-1)}
    \inf_{a>0}  a^{-\kappa+1} e^a 
    \left( 1 - a \int_0^\infty e^{-\lambda} (\lambda+a)^{-1}\,d\lambda \right)^{-1} \,.
\end{equation*}
Although, by inserting $K\leq (\kappa/S)^\kappa$, this yields a slightly worse bound on $L$ than the one in (\ref{eq:CLRconst}), in concrete applications one often has better bounds on $K$ available. In particular, in the case $T=-\Delta$ in $d=3$ one has $\kappa=3/2$ and $K=(4\pi)^{-3/2}$. The upper and lower bounds on $L$ derived this way then differ only by a factor $1.49$.
\end{remark}

Next, we turn to estimates on eigenvalue moments. Consider two sets of parameters $(\kappa,\gamma)\in (0,\infty)\times(0,\infty)$ and $(q,\theta)\in(2,\infty)\times(0,1)$ related by
\begin{equation}\label{eq:gammakappa}
 \gamma= \frac{q(1-\theta)}{q-2}\,, \qquad \kappa=\frac{q \theta}{q-2} \,,
\end{equation}
and
\begin{equation}\label{eq:qtheta}
 q = \frac{2(\gamma+\kappa)}{\gamma+\kappa-1} \,,
\qquad \theta = \frac{\kappa}{\gamma+\kappa} \,.
\end{equation}

\begin{theorem}[Equivalence of Sobolev and weak LT Inequalities]\label{ltwweighted}
Let $(\kappa,\gamma)$ and $(q,\theta)$ be as in \eqref{eq:gammakappa} and \eqref{eq:qtheta} and assume that $\gamma+\kappa>1$. Under Assumption~\ref{ass:bd} with $\kappa$ replaced by $\gamma+\kappa$ the following are equivalent:
\begin{enumerate}
 \item\label{it:gnweighted} $T$ satisfies a Sobolev interpolation inequality with exponent $q$, that is, there is a constant $S>0$ such that for all $u\in\mathrm{dom}\, t$,
\begin{equation}\label{eq:gnweighted}
t[u]^{\theta} \|u\|^{2(1-\theta)} \geq S \left( \int_X |u|^q \,dx \right)^{2/q}  \,.
\end{equation}
\item\label{it:ltwweighted} $T$ satisfies a weak LT inequality with exponent $\kappa$, that is, there is a constant $L>0$ such that for all $0\leq V\in L^{\gamma+\kappa}(X)$ and all $\tau>0$,
\begin{equation}\label{eq:wltweighted}
 N(-\tau,T-V) \leq L \ \tau^{-\gamma} \int_X V^{\gamma+\kappa} \,dx \,.
\end{equation}
\end{enumerate}
The respective sharp constants satisfy
\begin{equation}\label{eq:ltwconst}
(\theta^{-\theta} (1-\theta)^{-1+\theta} S)^{-\gamma-\kappa} \leq L \leq e^{\gamma+\kappa-1} (\theta^{-\theta} (1-\theta)^{-1+\theta} S)^{-\gamma-\kappa} \,.
\end{equation}
\end{theorem}

\begin{corollary}[LT Inequalities]\label{ltweighted}
 Let $T$ satisfy the Sobolev interpolation inequality \eqref{eq:gnweighted} for some $2<q<\infty$ and $0<\theta<1$ and let Assumption \ref{ass:bd} hold with $\kappa$ replaced by $q/(q-2)$. Define $0<\kappa<\infty$ and $0<\gamma<\infty$ by \eqref{eq:gammakappa}. Then for all $\tilde\gamma>\gamma$ and for all $0\leq V\in L^{\tilde\gamma+\kappa}(X)$ one has
\begin{equation}\label{eq:ltweighted}
 \mathrm{Tr}\,(T-V)_-^{\tilde\gamma} \leq L_{\tilde\gamma} \int_X V^{\tilde\gamma+\kappa} \,dx
\end{equation}
with
$$
L_{\tilde\gamma} \leq 
\frac{\tilde\gamma^{\tilde\gamma+1}}{\gamma^\gamma (\tilde\gamma-\gamma)^{\tilde\gamma-\gamma}}
\ \frac{\Gamma(\gamma+\kappa+1) \Gamma(\tilde\gamma-\gamma)}{\Gamma(\tilde\gamma+\kappa+1)} 
\ L
\,,
$$
where $L$ is the sharp constant in (\ref{eq:wltweighted}).
\end{corollary}


\subsection{Inclusion of magnetic fields}

The previous analysis can be extended to operators with magnetic fields, which do not satisfy Assumption~\ref{ass:bd}. A judicious use the diamagnetic inequality allows one to reduce the problem to the non-magnetic case. On the level of quadratic forms, the diamagnetic inequality means the following. Let $t$ and $t_A$ denote two closed, non-negative quadratic forms in $L^2(X)$. We say that $t_A$ satisfies a diamagnetic inequality with respect to $t$ if
for any $u\in\mathrm{dom}\, t_A$ and $v\in\mathrm{dom}\, t$ with $0\leq v\leq |u|$ one has $|u|\in\mathrm{dom}\, t$, $v\,\mathrm{sgn}\, u \in\mathrm{dom}\, t_A$ and 
\begin{equation}\label{eq:kato}
 t[v,|u|] \leq \text{Re}\ t_A[v \,\mathrm{sgn}\, u, u] \,.
\end{equation}
Here we use the definition $\mathrm{sgn}\, u(x):= u(x)/|u(x)|$ if $u(x)\neq 0$ and $\mathrm{sgn}\, u(x):=0$ if $u(x)=0$. Moreover, $t[\cdot,\cdot]$ denotes the sesqui-linear form associated to the quadratic form $t[\cdot]$ which is anti-linear in the first and linear in the second argument, and likewise for $t_A$. 

Let $T_A$ denote the operator corresponding to $t_A$.

\begin{theorem}[Inequalities with magnetic fields]\label{magnetic} Assume that $t$ satisfies Assumption~\ref{ass:bd} as well as either (\ref{eq:sobolevweighted}) or (\ref{eq:gnweighted}) and assume that $t_A$ satisfies a diamagnetic inequality with respect to $t$. Then the number of negative eigenvalues of $T_A-V$ satisfies the bounds (\ref{eq:clrweighted}) or (\ref{eq:wltweighted}) with the same upper bounds on the constants $L$ as in \eqref{eq:CLRconst} or \eqref{eq:ltwconst}.
\end{theorem}

Corollary~\ref{ltweighted} has a similar extension to the magnetic case as well. For further results about the magnetic version of CLR and LT inequalities we refer to \cite{Ro3, magest}.


\section{Illustrative examples}\label{sec:examples}

\subsection{The Laplacian}\label{sec:ltclass}

Let $T=-\Delta$ in $L^2(\mathbb R^d)$ corresponding to the quadratic form $t[u]:=\int_{\mathbb R^d} |\nabla u|^2\,dx$ with domain $H^1(\mathbb R^d)$. One easily checks that Assumption \ref{ass:bd} is satisfied for $\omega\equiv 1$. Moreover, if $2< q<\infty$ for $d=1,2$, or if $2<q\leq 2d/(d-2)$ for $d\geq 3$ the following Sobolev inequalities
\begin{equation}\label{eq:sobolev}
\left( \int_{\mathbb R^d} |\nabla u|^2\,dx \right)^{\theta} \left( \int_{\mathbb R^d} |u|^2\,dx \right)^{1-\theta}
\geq S_{q,d}  \left( \int_{\mathbb R^d} |u|^q\,dx \right)^{2/q}
\end{equation}
are well-known; see \cite[Sects.~8.3 and 8.5]{LiLo}. Here $\theta = d(\frac12-\frac 1q)$. By Theorem \ref{clrweighted} and Corollary~\ref{ltweighted} this implies the usual LT inequalities
\begin{equation}\label{eq:ltclass}
\mathrm{Tr}\left(-\Delta-V\right)_-^\gamma \leq L_{\gamma,d} \int_{\mathbb R^d} V_+^{\gamma+d/2} \,dx
\end{equation}
for $\gamma>(2-d)/2$ if $d=1,2$ and $\gamma\geq 0$ if $d\geq 3$. We
note that the inequality in the case $\gamma=1/2$, $d=1$, though being
valid \cite{HuLiTh,We}, cannot be obtained using the approach of
the present paper.  For a review of this topic, and remarks about the best constants,
see \cite{LaWe, Hu}.

If $A\in L^{2}_{\mathrm{loc}}(\mathbb R^d;\mathbb R^d)$ and $t_{A}[u] := \int_{\mathbb R^d} \left|(\nabla+iA) u\right|^2\,dx$, then \eqref{eq:kato} holds \cite{Ka,Si79}. Hence by Theorem \ref{magnetic}, inequality \eqref{eq:ltclass} holds with $-\Delta$ replaced by $-(\nabla+iA)^2$.


\subsection{Fractional Laplacians}

Let $f$ be a non-negative, differentiable function on $(0,\infty)$ such that $f'$ is completely monotone. We claim that the operator $T=f(-\Delta)$ satisfies Assumption \ref{ass:bd} with $\omega\equiv 1$. To verify this, we have to check by the Beurling--Deny theorem \cite[Sec. 1.3]{D} that $\exp(-tf(-\Delta))$ is positivity preserving and contractive in $L^\infty(\mathbb R^d)$. One easily checks that $E\mapsto e^{-tf(E)}$ is completely monotone and therefore by Bernstein's theorem \cite[Sec. I.5]{Do} it is of the form $e^{-tf(E)} = \int_0^\infty e^{-sE} \,d\mu_{f,t}(s)$ for some non-negative measure $\mu_{f,t}$. Hence $\exp(-tf(-\Delta)) = \int_0^\infty \exp(s\Delta) \,d\mu_{f,t}(s)$. Since $\exp(s\Delta)$ is positivity preserving and contractive in $L^\infty(\mathbb R^d)$ for any $s> 0$, so is $\exp(-tf(-\Delta))$. (For the contraction property we also use that $\int_0^\infty \,d\mu_{f,t}(s) = e^{-tf(0)}\leq 1$.)

In particular, the function $f(E)=E^s$ for $0<s<1$ is of the form described above. In this case, one has the Sobolev inequalities
\begin{align*}
 \left\|(-\Delta)^{s/2} u\right\|^{2\theta} \left(\int_{\mathbb R^d} | u|^2 \,dx \right)^{1-\theta}
 \geq S_{d,s,q} \left(\int_{\mathbb R^d} | u|^q \,dx \right)^{2/q}\,,
\qquad \theta = \frac d{2s}\left(\frac12-\frac1q\right) \,,
\end{align*}
for $2<q<\infty$ if $d\leq 2s$ and for $2<q\leq 2d/(d-2s)$ if $d>2s$. Hence Theorem \ref{clrweighted} and Corollary \ref{ltweighted} yield the inequalities
\begin{equation}
\mathrm{Tr}\left((-\Delta)^s-V\right)_-^\gamma \leq L_{\gamma,s,d} \int_{\mathbb R^d} V_+^{\gamma+d/2s} \,dx
\end{equation}
for $\gamma>(2s-d)/2s$ if $d\leq 2s$ and $\gamma\geq 0$ if $d>2s$. These inequalities appeared first in \cite{Da}.


\subsection{Periodic Schr\"odinger operators}

Let $W$ be a $\mathbb Z^d$-periodic function on $\mathbb R^d$ and consider the Schr\"odinger operator $-\Delta+W$ in $L^2(\mathbb R^d)$. With $E:=\inf\mathrm{spec}\,(-\Delta+W)$ the quadratic form $t[u] := \int_{\mathbb R^d} \left( |\nabla u|^2 + (W-E)|u|^2\right)\,dx$ is non-negative. Under very weak conditions on $W$ there is a periodic function $\omega$ satisfying $-\Delta \omega+W\omega=E\omega$, and $\omega$ is bounded (by elliptic regularity) and strictly positive (by Harnack's inequality). The representation
$$
t[u] = \int_{\mathbb R^d} |\nabla v|^2 \omega^2 \,dx\,,
\qquad u=\omega v\,,
$$
together with the properties of $\omega$ implies that Assumption \ref{ass:bd} is satisfied for any $\kappa>1$. Moreover, by the same representation the Sobolev inequalities \eqref{eq:sobolev} hold with the constant $S_{q,d}$ replaced by $(\inf\omega/\sup\omega)^2 S_{q,d}$. Therefore Theorem \ref{clrweighted} and Corollary~\ref{ltweighted} yield the CLR and LT inequalities
$$
\mathrm{Tr}\left(-\Delta +W-V-E\right)_-^\gamma \leq L_{\gamma,d}(W) \int_{\mathbb R^d} V_+^{\gamma+d/2} \,dx
$$
for the same values of $\gamma$ as in Subsection \ref{sec:ltclass}. This was first shown in \cite{ltcomp} by a different argument (which includes the case $\gamma=1/2$ and $d=1$).


\subsection{Hardy--Lieb--Thirring inequalities}

Let $d\geq 1$ and $0<s\leq 1$ such that $d>2s$, and denote by $\mathcal C_{s,d}$ the sharp constant in the Hardy inequality
$$
\int_{\mathbb R^d} |\xi|^{2s} |\hat u(\xi)|^2 \,d\xi \geq \mathcal C_{s,d} 
\int_{\mathbb R^d} |x|^{-2s} | u(x)|^2 \,dx \,,
\qquad u\in C_0^\infty(\mathbb R^d)\,.
$$
Here $\hat u(\xi) = (2\pi)^{-d/2} \int e^{-i\xi\cdot x} u(x)\,dx$ is the Fourier transform of $u$.
Explicitly (see \cite{He}), one has
$$
\mathcal C_{s,d} = 2^{2s} \frac{\Gamma((d+2s)/4)^2}{\Gamma((d-2s)/4)^2} \,.
$$
Let $t$ be the closure in $L^2(\mathbb R^d)$ of the non-negative quadratic form $\int_{\mathbb R^d} |\xi|^{2s} |\hat u(\xi)|^2 \,d\xi - \mathcal C_{s,d}  \int_{\mathbb R^d} |x|^{-2s} | u(x)|^2 \,dx$ defined on $C_0^\infty(\mathbb R^d)$. In \cite{ltpseudo} we have derived the following ground state representation formula,
\begin{align*}
t[u] = a_{s,d} \iint_{\mathbb R^d\times\mathbb R^d} \frac{|v(x)-v(y)|^2}{|x-y|^{d+2s}} \frac{dx}{ |x|^{(d-2s)/2}} \frac{dy}{ |y|^{(d-2s)/2}} \,, \qquad u=\omega v\,,
\end{align*}
where $\omega(x) = |x|^{-(d-2s)/2}$ and $a_{s,d}$ is a positive constant. This formula together with the fact that $C_0^\infty(\mathbb R^d\setminus\{0\})$ is a form core shows that Assumption \ref{ass:bd} holds. Moreover, in \cite{ltpseudo} we have shown that for any $2<q<2d/(d-2s)$ there is an $\tilde S_{d,s,q}>0$ such that
\begin{align*}
 t[u]^\theta \left(\int_{\mathbb R^d} | u|^2 \,dx \right)^{1-\theta}
 \geq \tilde S_{d,s,q} \left(\int_{\mathbb R^d} | u|^q \,dx \right)^{2/q}\,,
\qquad \theta = \frac d{2s}\left(\frac12-\frac1q\right) \,.
\end{align*}
In view of Corollary \ref{ltweighted} we obtain
\begin{equation}
 \label{eq:hlt}
\mathrm{Tr}\left((-\Delta)^s - \mathcal C_{s,d} |x|^{-2s} -V\right)_-^\gamma 
\leq \tilde L_{\gamma,d,s} \int_{\mathbb{R}^d} V_+^{\gamma+d/2s}\,dx
\end{equation}
for all $\gamma>0$ and the values of $s$ indicated above. This inequality for $s=1$ was first proved in \cite{liebhardy}. The proof sketched above is taken from \cite{ltpseudo}. For an alternative proof covering the cases $d\geq 3$ and $1<s<d/2$ we refer to \cite{hltsimple}.

Using Theorem \ref{magnetic} one can show that inequality \eqref{eq:hlt} holds also in the magnetic case, that is, with $(-\Delta)^s$ replaced by $|\nabla +iA|^{2s}$ for some $A\in L^{2}_{\mathrm{loc}}(\mathbb R^d;\mathbb R^d)$. This fact allowed us to prove stability of relativistic matter in magnetic fields up to the critical value of the nuclear charge; see \cite{ltpseudo} and also \cite{stability,LiSe}.


\section{Proofs of the main results}

\subsection{Proof of Theorem \ref{clrweighted} for $\omega\equiv 1$}

We begin by proving the easy implication $(\ref{it:clrweighted}\Rightarrow\ref{it:sobolevweighted})$. The CLR inequality and the variational principle imply that if $\int_X V^\kappa \,dx< L^{-1}$, then $t[u]\geq \int V|u|^2 \,dx$ for all $u\in\mathrm{dom}\, t$. Choosing $V= \alpha_\varepsilon |u|^{2/(\kappa-1)}$ with
$$
\alpha_\varepsilon =(1-\varepsilon) \left( L \int_X |u|^{2\kappa/(\kappa-1)} \,dx \right)^{-1/\kappa}
$$
and letting $\varepsilon\to 0$ we obtain \eqref{eq:sobolevweighted} with $S\geq L^{-1/\kappa}$.

\bigskip
Next, we shall present two proofs of the implication $(\ref{it:sobolevweighted}\Rightarrow\ref{it:clrweighted})$. The first one is an abstraction of the semi-group proof of the CLR inequality in \cite{Li}. It relies on the heat kernel bound
\begin{equation}\label{eq:heatkernel}
 \|\exp(-sT)\|_{1\to\infty} \leq K s^{-\kappa} \,.
\end{equation}
We recall that by Varopoulos' theorem \cite{Va} \eqref{eq:sobolevweighted} is equivalent to \eqref{eq:heatkernel}. An abstract version of an argument by Nash \cite{Na} allows us to derive \eqref{eq:heatkernel} from \eqref{eq:sobolevweighted} with constant
\begin{equation}\label{eq:heatconst}
K \leq (\kappa/S)^\kappa \,.
\end{equation}
Indeed, the Sobolev inequality \eqref{eq:sobolevweighted} implies via H\"older the Nash inequality
\begin{equation}
 t[u]^{\frac{q}{2(q-1)}} \left(\int_X |u| \,dx \right)^{\frac{q-2}{q-1}}
\geq S^{\frac{q}{2(q-1)}} \int_X |u|^2 \,dx \,.
\end{equation}
By Nash's argument (see \cite[Sec.~2.4]{D} or \cite[Thm.~8.16]{LiLo}), using the contraction property in $L^1$, this implies
\begin{equation}
 \|\exp(-sT)\|_{1\to2}^2 \leq \left( \kappa / 2S \right)^\kappa s^{-\kappa} \,,
\end{equation}
which yields \eqref{eq:heatkernel} and \eqref{eq:heatconst} by duality and the semi-group property.

With \eqref{eq:heatkernel} at hand we can now follow the arguments in \cite{Li}, replacing path integrals by Trotter's product formula (see also \cite{RoSo}). Defining for any non-negative, lower semi-continuous function $f$ on $\mathbb{R}_+$ with $f(0)=0$
\begin{equation}\label{eq:f}
      F(\lambda):=\int_0^\infty f(\mu) e^{-\mu/\lambda}\mu^{-1}\,d\mu \,,
      \qquad \lambda >0 \,,
\end{equation}
one has the trace formula
    \begin{equation}\label{eq:traceformula}
       \begin{split}
         & \mathrm{Tr}\,  F(V^{1/2}\, T^{-1} V^{1/2})
         = \int_0^\infty \frac{ds}s  \lim_{n\to\infty}
         \int_X \cdots\int_X \, dx_1 \cdots dx_n \\
         & \qquad\qquad\qquad\qquad\qquad  \times \prod_{j=1}^n k\left(x_j,x_
{j-1},\frac sn\right)
         f\left(\frac sn \sum_{k=1}^n V(x_k) \right)
       \end{split}
     \end{equation}
with the convention that $x_0=x_n$ and $k(x,y,s)=\exp(-sT)(x,y)$. By Assumption~\ref{ass:bd} and the Beurling--Deny theorem (see, e.g., \cite[Sec.~1.3]{D}) $k$ is non-negative. If, in addition, $f$ is convex then we can bound 
$f\left(\frac sn \sum V(x_k) \right) \leq  \frac1n\sum f(sV(x_k))$ and obtain, using the semi-group property,
$$
\mathrm{Tr}\,  F(V^{1/2}\, T^{-1} V^{1/2})
         \leq \int_0^\infty \frac{ds}s \int_X dx\, k(x,x,s) f(sV(x)) \,.
$$
Now the heat kernel decay \eqref{eq:heatkernel} implies
$$
\mathrm{Tr}\,  F(V^{1/2}\, T^{-1} V^{1/2})
\leq K \int_0^\infty \frac{ds}{s^{\kappa+1}} \int_X dx f(sV(x))
= K \int_X V^\kappa \,dx  \int_0^\infty f(\mu) \frac{d\mu}{\mu^{\kappa+1}} \,.
$$
By the Birman-Schwinger principle, $N(0,T-V)$ coincides with the number of eigenvalues larger than one of the operator $V^{1/2}\, T^{-1} V^{1/2}$. Hence, since $F$ is increasing,
$$
N(0,T-V) \leq F(1)^{-1} \mathrm{Tr}\,  F(V^{1/2}\, T^{-1} V^{1/2}) 
\leq K \int_X V^\kappa \,dx \, F(1)^{-1} \int_0^\infty f(\mu) \frac{d\mu}{\mu^{\kappa+1}} \,,
$$
and the sought bound follows by choosing $f(\mu)=(\mu-a)_+$ and optimizing over $a>0$.

The only place where part (3) of Assumption~\ref{ass:bd} entered in the proof is to obtain the heat kernel bound (\ref{eq:heatkernel}) from the Sobolev inequality (\ref{eq:sobolevweighted}). This part of the assumption can thus be omitted if one is able to obtain such a bound by other means.

\bigskip
The second proof of the implication $(\ref{it:sobolevweighted}\Rightarrow\ref{it:clrweighted})$ is an abstraction of Li and Yau's proof \cite{LYa} of the CLR inequality and its improvement in \cite{BlReSt}. By an approximation argument we may assume that $V\in L^1\cap L^\infty$ and $V>0$ a.e. Moreover, for the sake of simplicity we assume that the embedding of the completion of $\mathrm{dom}\, t$ with respect to $t$ into $L^q(X)$ is injective; see \cite{LeSo} for an additional argument in the general case. We consider the non-negative operator $\Upsilon$ in $L^2(X,V)$ given by the quadratic form $t[v]$. We shall prove that
\begin{equation}\label{eq:liyau}
 \mathrm{Tr}\, (2\Upsilon)^{-1}\exp(-2s\Upsilon) \leq (\kappa-1)^{\kappa-1} (2S)^{-\kappa} \int_X V^\kappa \,dx \ s^{-\kappa+1} \,.
\end{equation}
Since $\Upsilon$ in $L^2(X,V)$ is unitarily equivalent to the inverse of the Birman--Schwinger operator $V^{1/2}T^{-1}V^{1/2}$ in $L^2(X)$ one has $N(T-V)=N(1,\Upsilon^{-1})$. The inequality $N(1,\Upsilon^{-1}) \leq 2e^{2t} \mathrm{Tr}\, (2\Upsilon)^{-1}\exp(-2t\Upsilon)$, together with \eqref{eq:liyau} and optimization in $t$, will then imply \eqref{eq:clrweighted} with the upper bound on $L$ stated in \eqref{eq:CLRconst}.

In order to prove \eqref{eq:liyau} we consider the operators $H_\beta(s) = (2\Upsilon)^{-\beta} \exp(-s\Upsilon)$ for $\beta\geq 0$. From the Sobolev inequality (\ref{eq:sobolevweighted}) and Assumption~\ref{ass:bd} one concludes, as in \cite{LeSo}, that $H_0(s)$, and hence also $H_\beta(s)$, are integral operators with non-negative kernels $H_\beta(x,y,s)$. We abbreviate $h_\beta(s):=\mathrm{Tr}\, H_\beta(2s)$ and estimate, using H\"older with $\frac1q+\frac{q-2}q+\frac1q=1$,
\begin{align*}
h_\beta(s)&=\int_X dx\,V(x) \left( \int_X dy\, V(y) H_\beta(y,x,s) H_0(x,y,s) \right) \\
&\leq \int_X dx\,V(x) \left( \int_X dy\, H_\beta(y,x,s)^q \right)^{\frac 1q}
\left( \int_X dy\, H_0(x,y,s) V(y)^{\frac{q-1}{q-2}} \right)^{\frac{q-2}q} \\
& \qquad \qquad\qquad\quad \times \left( \int_X dy\, H_0(x,y,s)^2 V(y) \right)^{\frac 1q} \,.
\end{align*}
Using H\"older once more with $\frac12+\frac{q-2}{2q}+\frac1q=1$ we obtain
\begin{align*}
h_\beta(t)
\leq & \left( \int_X dx\,V(x) \left( \int_X dy\, H_\beta(y,x,s)^q \right)^{\frac 2q} \right)^{\frac12} \\
& \times \left( \int_X dx\,V(x) \left( \int_X dy\, H_0(x,y,s) V(y)^{\frac{q-1}{q-2}} \right)^2 \right)^{\frac{q-2}{2q}} \\
& \times \left( \int_X dx\,V(x) \int_X dy\, H_0(x,y,s)^2 V(y) \right)^{\frac1q}
=: A^{\frac12} B^{\frac{q-2}{2q}} C^{\frac1q} \,.
\end{align*}
We estimate $A$ by the Sobolev inequality \eqref{eq:sobolevweighted},
$$
A \leq S^{-1} \int_X dx\,V(x) \, t[H_\beta(\cdot,x,s)] = (2S)^{-1} \, h_{2\beta-1}(s) \,.
$$
The contraction property of $\exp(-s\Upsilon)$ in $L^2(X,V)$ implies
$$
B = \| \exp(-s\Upsilon) V^{\frac{1}{q-2}} \|^2_{L^2(X,V)} \leq \| V^{\frac{1}{q-2}} \|^2_{L^2(X,V)} = \int_X V^\kappa \,dx \,.
$$
Moreover, $C=h_0(s)$. Hence, choosing $\beta=1$ and using $h_0(s)=-h_1'(s)$, we have shown
$$
h_1(s) \leq (2S)^{-1} \left(\int_X V^\kappa \,dx\right)^{\frac 1\kappa} \left(-h_1'(s)\right)^{\frac 2q} \,
$$
which implies \eqref{eq:liyau} and completes the sketch of the proof.

Note that the only place where part (3) of Assumption~\ref{ass:bd} entered in the second proof is the existence of integral kernels for the operators $H_\beta(s)$. Hence this part of the  assumption can, in principle, be omitted if this property can be shown by other means.


\subsection{Proof of Theorem \ref{clrweighted} for arbitrary $\omega$}

Since the proof of the implication $(\ref{it:clrweighted}\Rightarrow\ref{it:sobolevweighted})$ in the previous subsection did not use Assumption \ref{ass:bd} we are left with proving $(\ref{it:sobolevweighted}\Rightarrow\ref{it:clrweighted})$. We will deduce this from the case $\omega\equiv 1$.
We may assume that $T$ is positive definite for otherwise we consider $T+\varepsilon$ and let $\varepsilon\to 0$ in the inequality obtained. The quadratic form $t_\omega[v] := t[\omega v]$ with domain $\omega^{-1} \mathcal D$ is closable
in the Hilbert space $L^2(X,d\mu)$ with measure $d\mu := \omega^{2\kappa/(\kappa-1)} dx$. (Here we use that $T$ is positive definite and that $t$ is closed.) Let $T_\omega$ be the corresponding self-adjoint operator in $L^2(X,d\mu)$. We note that $t_\omega$ satisfies Assumption \ref{ass:bd} with $\omega\equiv 1$ (it suffices to verify this assumption on a form core, see \cite[Lem. 1.3.4]{D}) and that the Sobolev inequality \eqref{eq:sobolevweighted} for $t$ can be written as
$$
t_\omega[v] \geq S \left(\int_X |v|^{q} \,d\mu \right)^{2/q} \,.
$$
Moreover, by the variational principle,
\begin{align*}
 N(T-V) & = \sup\left\{\dim M : \ M\subset\mathcal D, \, t[u] < \int_X V |u|^2 \,dx \ \text{for all}\ 0\not\equiv u\in M \right\} \\
& = \sup\left\{\dim \tilde M : \ \tilde M\subset\omega^{-1} \mathcal D, \, t_\omega[v] < \int_X \tilde V |v|^2 \,d\mu \ \text{for all}\ 0\not\equiv v\in \tilde M \right\} \\
& = N(T_\omega - \tilde V)
\end{align*}
where $\tilde V := \omega^{-2/(\kappa-1)} V$. Since
$$
\int_X \tilde V^{\kappa} \,d\mu = \int_X V^{\kappa} \,dx \,,
$$
the assertion follows from the $\omega\equiv 1$ case of Theorem \ref{clrweighted}.


\subsection{Proof of Theorem \ref{ltwweighted}}\label{sec:omega}

We shall deduce the result of Theorem \ref{ltwweighted} for positive $\gamma$ from that of Theorem \ref{clrweighted} for $\gamma=0$. To do so, we consider the operator $T_\tau := \tau^{-1+\theta}(T+\tau)$ and its quadratic form $t_\tau$. Then condition \eqref{eq:wltweighted} is equivalent to
\begin{equation*}
N(T_\tau - V) \leq L \int_X V^{\gamma+\kappa} \,dx \,,
\qquad \tau>0 \,,
\end{equation*}
for all $0\leq V\in L^{\gamma+\kappa}(X)$. Moreover, using that for $\alpha,\beta>0$
\begin{equation*}
\min_{\tau>0} (\alpha \tau^{-1+\theta} + \beta \tau^\theta) 
= \theta^{-\theta} (1-\theta)^{-1+\theta} \alpha^{\theta}\beta^{1-\theta} \,,
\end{equation*}
condition \eqref{eq:gnweighted} is equivalent to
\begin{equation*}
t_\tau[u] \geq \theta^{-\theta} (1-\theta)^{-1+\theta} S \left( \int_X |u|^q \,dx \right)^{2/q} \,,
\qquad \tau>0 \,.
\end{equation*}
Noting that $T_\tau$ satisfies Assumption \ref{ass:bd}, the assertion follows from Theorem \ref{clrweighted}.


\subsection{Proof of Corollary \ref{ltweighted}}

By Theorem \ref{ltwweighted} the Sobolev interpolation inequality \eqref{eq:gnweighted} implies the weak LT inequality \eqref{eq:wltweighted}. We shall now use an interpolation argument from \cite{LiTh} in order to deduce the strong LT inequality for $\tilde\gamma>\gamma$ from a weak LT inequality for $\gamma$. For any fixed $0<s<1$ the variational principle implies
\begin{equation*}
N(-\tau,T-V) \leq N(-(1-s)\tau, T-(V-s\tau)_+)\,.
\end{equation*}
Hence the representation
\begin{equation*}\label{eq:momentintegral}
\mathrm{Tr}\, (T-V)_-^{\tilde\gamma} = \tilde\gamma \int_0^\infty N(-\tau,T-V) \tau^{\tilde\gamma-1}\,d\tau
\end{equation*}
together with the weak LT inequality implies that
\begin{align*}
 \mathrm{Tr}\,(T-V)_-^{\tilde\gamma} 
& \leq L \tilde\gamma (1-s)^{-\gamma} \int_0^\infty \int_X (V-s\tau)_+^{\gamma+\kappa}\,dx \tau^{\tilde\gamma-\gamma-1}\,d\tau \\
& = L \tilde\gamma (1-s)^{-\gamma}  s^{-\tilde\gamma+\gamma} B(\gamma+\kappa+1,\tilde\gamma-\gamma) \int_X V^{\tilde\gamma+\kappa}\,dx \,,
\end{align*}
with $B(\cdot,\cdot)$ the beta function. Minimizing in $s\in (0,1)$ yields the claimed inequality with constant
$$
L_{\tilde\gamma} \leq 
\frac{\tilde\gamma^{\tilde\gamma+1}}{\gamma^\gamma (\tilde\gamma-\gamma)^{\tilde\gamma-\gamma}}
\ \frac{\Gamma(\gamma+\kappa+1) \Gamma(\tilde\gamma-\gamma)}{\Gamma(\tilde\gamma+\kappa+1)} 
\ L
\,.
$$


\subsection{Proof of Theorem~\ref{magnetic}}

{F}rom the diamagnetic inequality (\ref{eq:kato}) for the quadratic forms one concludes that the heat kernel for the operator $T_A$ is, in absolute value, pointwise bounded from above by the heat kernel for $T$. This was proved by Hess, Schrader, Uhlenbrock \cite{HeScUh} and Simon \cite{Si}; see \cite[Sec. 2]{Ou} for a quadratic form version of this result. For $\omega\equiv 1$ this immediately implies that the first proof of $(\ref{it:sobolevweighted}\Rightarrow\ref{it:clrweighted})$ in Theorem \ref{clrweighted}, using the method in \cite{Li}, extends to the magnetic case with the same bound on the constant. For general $\omega$ one proceeds as in Subsection \ref{sec:omega}, noting that $t_{A,\omega}$ satisfies a diamagnetic inequality in the sense of \eqref{eq:kato} with respect to $t_\omega$.

A similar argument shows that the operator $\Upsilon_A$ in the second proof of $(\ref{it:sobolevweighted}\Rightarrow\ref{it:clrweighted})$ in Theorem \ref{clrweighted}, using the method in \cite{LYa}, satisfies a diamagnetic inequality with respect to $\Upsilon$. Hence
$$
\mathrm{Tr}\, (2\Upsilon_A)^{-1}\exp(-2t\Upsilon_A) \leq \mathrm{Tr}\, (2\Upsilon)^{-1}\exp(-2t\Upsilon) \,.
$$
Hence \eqref{eq:liyau} leads to the same estimate in the magnetic case as in the non-magnetic case.


\bibliographystyle{amsalpha}

\end{document}